\documentclass[letterpaper,11pt]{amsart}

\usepackage[all]{xy}                        %

\CompileMatrices                            

\UseTips                                    

\input xypic
\usepackage[bookmarks=true]{hyperref}       

\usepackage{amssymb,latexsym,amsmath,amscd}
\usepackage{xspace}

\usepackage{graphicx}

\reversemarginpar

\theoremstyle{plain}
\newtheorem{theorem}{Theorem}[section]
\newtheorem*{theorem*}{Theorem}
\newtheorem{proposition}[theorem]{Proposition}

\newtheorem{lemma}[theorem]{Lemma}

\theoremstyle{definition}

\newtheorem{remark}[theorem]{Remark}

\newcommand{\enm}[1]{\ensuremath{#1}}          %
\newcommand{\op}[1]{\operatorname{#1}}
\newcommand{\cal}[1]{\mathcal{#1}}

\newcommand{\CC}{\enm{\mathbb{C}}}

\newcommand{\FF}{\enm{\mathbb{F}}}

\newcommand{\PP}{\enm{\mathbb{P}}}

\newcommand{\Dd}{\enm{\cal{D}}}
\newcommand{\Ee}{\enm{\cal{E}}}

\newcommand{\Oo}{\enm{\cal{O}}}
\newcommand{\Pp}{\enm{\cal{P}}}

\newcommand{\Ww}{\enm{\cal{W}}}

\renewcommand{\phi}{\varphi}
\renewcommand{\theta}{\vartheta}
\renewcommand{\epsilon}{\varepsilon}

\newcommand{\Pic}{\op{Pic}}

\newcommand{\Hom}{\op{Hom}}
\newcommand{\Ext}{\op{Ext}}
\newcommand{\Sec}{\op{Sec}}

\newcommand{\old}[1]{}

\begin{document}

\title[A Note on a Brill-Noether Locus]{A Note on a Brill-Noether Locus\\ over a non-hyperelliptic curve of genus 4}
\author{Sukmoon Huh}
\address{Korea Institute for Advanced Study \\
Hoegiro 87, Dongdaemun-gu \\
Seoul 130-722, Korea}
\email{sukmoon.huh@math.unizh.ch}
\thanks{The author would like to thank Edoardo Ballico and Robert Pignatelli for many advices.}
\begin{abstract}
We prove that a certain Brill-Noether locus over a non-hyperelliptic curve $C$ of genus 4, is isomorphic to the \textit{Donagi-Izadi cubic threefold} in the case when the pencils of the two trigonal line bundles of $C$ coincide.
\end{abstract}

\maketitle

\section{Introduction}
Let $C$ be a non-hyperelliptic curve of genus 4 over $\CC$ and then $C$ is embedded into $\PP_3$ by the canonical embedding and there exists a unique quadric surface $Q\subset \PP_3$ containing $C$. If we let $g_3^1$ and $h_3^1$ be the two trigonal line bundles such that $g_3^1\otimes h_3^1 = \Oo_C(K_C)$, the canonical line bundle, then $Q$ is singular if and only if the two pencils $|g_3^1|$ and $|h_3^1|$ coincide.

Let $SU_C(2, K_C)$ be the moduli space of semi-stable bundles of rank 2 on $C$ with the canonical determinant and $\Ww^r$ be the Brill-Noether locus defined as the closure of the set of stable bundles $E$ with $h^0(E)\geq r+1$. In \cite{OPP}, $\Ww^2$ was proven to be isomorphic to the \textit{Donagi-Izadi cubic threefold}. In \cite{Huh4}, we gave a different proof of this when $Q$ is smooth, using the fact that the moduli space of stable sheaves of rank 2 on $Q$ with the Chern classes $c_1=\Oo_Q(1,1)$ and $c_2=2$, is isomorphic to $\PP_3$.

In this article, we use the same trick of \cite{Huh4} to the Hirzebruch surface $\FF_2$ and derive the same result on $\Ww^2$ when $Q$ is a quadric cone in $\PP_3$.  Unlike the situation in \cite{Huh4}, the determinant of sheaves that we choose over $\FF_2$ is not ample, which prevents us from using the definition of stability. Instead, we use a parametrization $\Pp$ of the vector bundles on $\FF_2$ admitting a certain exact sequence. We show that $\Pp$ is isomorphic to $\PP_3$, the original ambient space into which $C$ is embedded by the canonical embedding. From the investigation of this parametrization, we show that the restriction map from $\Pp$ to $\Ww^2$ is given by the complete linear system $|I_C(3)|$, implying that $\Ww^2$ is isomorphic to the \textit{Donagi-Izadi cubic threefold}.

\section{Main Theorem}
Let $\FF_2=\PP ( F)$ be the Hirzebruch surface with a section $\sigma$ whose self-intersection is $-2$, where $F\simeq \Oo_{\PP_1}\oplus \Oo_{\PP_1}(2)$. Recall that $\FF_2$ is the minimal resolution of the quadric cone $Q\subset \PP_3$ at the vertex point $P_0$. The section $\sigma$ is the exceptional curve of the resolution. Let $f$ be a fibre of the ruling $\pi : \FF_2 \rightarrow \PP_1$ and then $\Pic (\FF_2)$ is freely generated by $\sigma$ and $f$. We will denote the line bundle $\Oo_{\FF_2}(a\sigma + bf)$ by $\Oo(a,b)$ and $E\otimes \Oo(a,b)$ by $E(a,b)$ for a coherent sheaf $E$ on $\FF_2$. Note that the canonical line bundle is $\Oo(-2,-4)$. Then the resolution $\phi : \FF_2 \rightarrow Q$ is given by the linear system $|\Oo(1,2)|$. Here, the line bundle $H=\Oo(1,2)$ is the tautological line bundle $\Oo_{\FF_2}(1)$ on $\FF_2$, which is nef but not ample.

\begin{lemma}\label{coh}
We have
$$H^i(\Oo_{\FF_2}(aH+bf))=\left\{
                            \begin{array}{ll}
                              0, & \hbox{if $a=-1$;} \\
                              H^i(\PP_1, S^a(F)\otimes \Oo_{\PP_1}(b)), & \hbox{if $a\geq 0$;} \\
                              H^{2-i}(\PP_1, S^{-2-a}(F)\otimes \Oo_{\PP_1}(-b)), & \hbox{if $a\leq -2$ .}
                            \end{array}
                          \right.$$
\end{lemma}
\begin{proof}
From the Leray spectral sequence,
$$H^i(\PP_1, R^j\pi_*\Oo_{\FF_2}(aH+bf))\Rightarrow H^{i+j}(\FF_2, \Oo_{\FF_2}(aH+bf))$$
and $R^i\pi_*\Oo_{\FF_2}(aH+bf)=R^i\pi_*\Oo_{\FF_2}(aH)\otimes \Oo_{\PP_1}(b)=0$ for $i>0$ and $a\geq -1$, we have
$$H^i(\FF_2, \Oo_{\FF_2}(aH+bf))=H^i(\PP_1, \pi_*\Oo_{\FF_2}(aH+bf)),$$
for all $a\geq -1$. Since $\pi_*\Oo_{\FF_2}(aH+bf)=S^a(F)\otimes \Oo_{\PP_1}(b)$ if $a\geq 0$ and 0 otherwise, we get the first and second assertions. The last case can be derived from the second case, using the Serre duality.
\end{proof}

Let $\Pp$ be the set of non-trivial sheaves of rank 2 on $\FF_2$ with the Chern classes $c_1=\Oo(1,2)$ and $c_2=2$, which are fitted into the following exact sequence,
\begin{equation}\label{ext1}
0\rightarrow \Oo(1,0) \rightarrow E \rightarrow \Oo(0,2)\rightarrow 0.
\end{equation}
 Note that $h^0(E)$ is 3 or 4 from the sequence (\ref{ext1}). If we let $\Oo(a,b)$ be a sub-bundle of $E$ from a section $s\in H^0(E)$, we have $(a,b)=(0,0),(0,1),(0,2),(1,0)$ or $(1,1)$, since $c_2(E)=2$. In the case of $(a,b)=(0,2)$, $E$ turns out to be isomorphic to $\Oo(1,0)\oplus \Oo(0,2)$, which is excluded.
In the case of $(1,1)$, $E$ is fitted into
\begin{equation}\label{ext2}
0\rightarrow \Oo(1,1) \rightarrow E \rightarrow I_p(0,1)\rightarrow 0,
\end{equation}
where $p$ is a point on $\FF_2$, implying that $h^0(E(-1,0))=2$. But this is not true since $h^0(E(-1,0))=1$ from (\ref{ext1}). Thus we obtain only $\Oo(0,0)$, $\Oo(1,0)$ or $\Oo(0,1)$ as sub-bundles of $E\in \Pp$ from sections of $E$.

Let $\Dd$ be the set of sheaves $E\in \Pp$, fitted into the following exact sequence,
\begin{equation}\label{ext3}
0\rightarrow \Oo(0,1) \rightarrow E \rightarrow I_p(1,1) \rightarrow 0,
\end{equation}
where $p$ is a point on $\FF_2$. Since the dimension of $\Ext^1 (I_p(1,1),\Oo(0,1))$ is 1, so we have the unique non-trivial extension of (\ref{ext3}) to each point $p\in \FF_2$ (It can be easily checked that the trivial extension does not lie in $\Pp$). Since $h^0(E(-1,0))=1$, we have $E\in \Pp$ and thus $\Dd \subset \Pp$. We also have $h^0(E(0,-1))=1+h^0(I_p(1,0))=2$ if $p  \in \sigma$ and $1$ otherwise. Similarly, $h^0(E)=4$ if $p \in \sigma$ and $3$ otherwise.

Let $(a,b)=(0,0)$ and so we have the exact sequence,
\begin{equation}\label{ext4}
0\rightarrow \Oo \rightarrow E \rightarrow I_Z(1,2) \rightarrow 0,
\end{equation}
where $Z$ is a 0-cycle on $\FF_2$ with length 2. Let us denote the extension classes of type (\ref{ext4}) by
$$\PP(Z):= \PP \Ext^1(I_Z(1,2), \Oo),$$
then $\PP(Z)$ is isomorphic to $\PP H^0(\Oo_Z)^*\simeq \PP_1$. Since $1=H^0(E(-1,0))=H^0(I_Z(0,2))$, we have two fibres $f_1$ and $f_2$ of $\pi$, each containing a point of $Z$. In fact, in the case when $Z$ is contained in a fibre $f$, we have the extension (\ref{ext2}).

\begin{proposition}
We have the following descriptions on $\Pp$:
\begin{enumerate}
\item $\Pp$ is isomorphic to $\PP_3$.
\item $\Dd\subset \Pp$ is a quadric cone $Q'$
\item The vertex $P_0'$ of $\Dd$ corresponds to the unique vector bundle $E_0\in \Pp$ such that $h^0(E_0)=4$
\end{enumerate}
\end{proposition}
\begin{proof}
The assertion (1) is clear since $h^0(E(-1,0))=1$ for all $E\in \Pp$ and $\PP \Ext^1 (\Oo(0,2), \Oo(1,0))$ is isomorphic to $\PP H^1(\Oo(1,-2))\simeq \PP_3$.
Now there exists a universal extension
$$0\rightarrow q^* \Oo(1,0) \rightarrow \Ee \rightarrow q^*\Oo(0,2) \rightarrow 0,$$
on $\Pp \times \FF_2$ ($q$ is the projection to $\FF_2$) such that $\Ee|_{ \{p\}\times \FF_2}$ is isomorphic to an extension corresponding to $p\in \Pp$. Let $\Ee '$ be an extension of $I_{\triangle} \otimes q^*\Oo(1,1)$ by $q^* \Oo(0,1)$ over $\FF_2 \times \FF_2$ such that the restriction of $\Ee '$ to $\{p\}\times \FF_2$ is the unique non-trivial extension of $I_p(1,1)$ by $\Oo(0,1)$. Here, $\triangle$ is the diagonal of $\FF_2 \times \FF_2$. The existence of such $\Ee '$ is guaranteed because
\begin{align*}
   & H^3(I_{\triangle}\otimes p^*\Oo(-1,-4) \otimes q^*\Oo(-2,-4)) \\ \simeq& H^2(\Oo_{\triangle} \otimes p^*\Oo(-1,-4) \otimes q^*\Oo(-2,-4))\\
 \simeq& H^2(\FF_2, \Oo(-1,-4)\otimes p_*q^*\Oo(-2,-4))\\
 \simeq& H^2(\Oo(-3,-8)) \simeq H^0(\Oo(1,4))
\end{align*}
is not zero. Since each restriction to $\{p\}\times \FF_2$ is contained in $\Pp$, we have a morphism $\chi$ from $\FF_2$ to $\Pp$ and the image of $\chi$ is $\Dd$. Now assume that $h^0(E)=4$ for $E\in \Pp$. It can be easily checked that there exists a section of $E$ for which $E$ is fitted into (\ref{ext4}). Thus, $4=h^0(E)=1+h^0(I_Z(1,2))$, i.e. $h^0(I_Z(1,2))=3$. This implies that $Z$ is contained in $\sigma$. From (\ref{ext4}), we also have $h^0(E(0,-1))>0$. In particular, $E$ is also fitted into (\ref{ext3}) with $p\in \sigma$.

Let $s_1, s_2$ be two sections of $E(-1,0)$ such that $p_1$ is the only zero of $s_1$ and $s_2$. If $s_1$ and $s_2$ are different, we can find $p_2\not= p_1$ such that $as_1+bs_2$ is zero at $p_2$ for some $a,b\not=0$, which is absurd because $p_1$ is also the unique zero of $as_1+bs_2$. Thus for all $p_1\in \sigma$, we have the unique $E$ such that $h^0(E)=4$. In particular, the map $\chi$ contracts $\sigma$ to a point in $\Pp$. Let $E\in \Dd$. If $p \not \in \sigma$, we have $h^0(E(0,-1))=1$ so that we can assign a different $E$ for each $p\not \in \sigma$. Thus $\chi$ is the minimal resolution of a quadric cone $Q'\subset \Pp$ at the vertex point $P_0'$ corresponding to the sheaf $E_0$ admitting (\ref{ext3}) with $p\in \sigma$.
\end{proof}

\begin{remark}\label{rmk}{$ $}
\begin{enumerate}
\item Let us consider the definition of stability on the sheaves of rank 2 on $\FF_2$ with the Chern classes $c_1=\Oo(1,2)$ and $c_2=2$ with respect to the nef divisor $H=\Oo(1,2)$. It can be checked that such sheaves admits an exact sequence (\ref{ext1}). Since all the sheaves in $\Dd$ contains $\Oo(0,1)$ as sub-bundle, it contradicts to the stability condition. So the space of stable sheaves in this sense, is isomorphic to $\PP_3 \backslash Q$ and in particular, it is not projective.
\item Let us assume that a non-trivial bundle $E$ with the extension (\ref{ext3}), $p\in \FF_2 \backslash \sigma$, admits an extension (\ref{ext4}) with $Z\in \FF_2^{[2]}$, where $\FF_2^{[2]}$ be the Hilbert scheme of 0-cycles of length 2 on $\FF_2$. In these two extensions, $\Oo$ is a sub-bundle of $\Oo(0,1)$, otherwise, $E$ containes $\Oo\oplus \Oo(0,1)$ as a sub-bundle, which is absurd. Thus we have a surjection from $I_Z(1,2)$ to $I_p(1,1)$. In particular, $\Hom (I_Z(1,2), I_p(1,1))$ is non-trivial. As a result, if we take $\Hom (\cdot , I_p(1,1))$ to the exact sequence,
$$0\rightarrow I_Z(1,2) \rightarrow \Oo(1,2) \rightarrow \Oo_Z \rightarrow 0,$$
then we know that $\Ext^1 (\Oo_Z, I_p)$ is non-trivial, which implies that $p\in Z$. Let us denote this $p$ by $p_E$.
\end{enumerate}
\end{remark}

For $E\in \Pp$, we consider the determinant map
$$\lambda_E : \wedge^2 H^0(E) \rightarrow H^0(\Oo(1,2)).$$
Recall that the dimension of $\wedge^2 H^0(E)$ is 3 if $E \not= E_0$.

\begin{lemma}\label{det}
If $E \in \Pp \backslash \Dd$, then $\lambda_E$ is injective.
\end{lemma}
\begin{proof}
Let $s_1$ and $s_2$ be two sections of $H^0(E)$ for which $s_1 \wedge s_2$ is a non-trivial element in $\ker (\lambda_E)$. It would generate a subsheaf $F$ of $E$ such that $h^0(F)\geq 2$. Since $c_2(E)=2$, it can be easily checked that the only possibility for $F$ is $\Oo(0,1)$ or $I_p(1,1)$, where $p$ is a point on $\FF_2$ with the following exact sequence,
$$0\rightarrow I_p(1,1) \rightarrow E \rightarrow \Oo(0,1) \rightarrow 0.$$
Let us assume that $E \not \in \Dd$ and in particular, $I_p(1,1)$ is the only possibility for $F$. From the previous result, $E$ is locally free. Since $\Oo(0,1)$ is torsion-free, so $I_p(1,1)$ must be a line bundle, which is absurd.
\end{proof}

Let us denote by $p_E\in \PP_3$, the point corresponding to the dual of the cokernel of $\lambda_E$. As vector subspaces of $H^0(\Oo(1,2))$, we see that $H^0(I_Z(1,2))$ is contained in the image of $\wedge^2 H^0(E)$ and so $p_E$ is contained in $H^0(\Oo_Z)^*$ as a vector subspace of $H^0(\Oo(1,2))^*$. It implies that $p_E$ is a point in $\PP_3$, contained in all secant lines of $\phi(Z)$ for which $E$ admits an extension (\ref{ext4}). This argument with the remark \ref{rmk}, gives us a map from $\eta : \Pp \rightarrow \PP_3$ sending $E$ to $p_E$ for $E\in \Pp \backslash \{E_0\}$. Clearly, this map extends to $E_0$ by assigning the vertex $P_0 \in Q$ because in the extension (\ref{ext4}) of $E_0$, the support of $Z$ should lie on $\sigma$ due to the fact that $h^0(E)=4$ and so $h^0(I_Z(1,1))=2$. Note that, for $E\in \Pp \backslash \Dd$, $p_E$ lies outside $Q$ and so we get the following statement.

\begin{proposition}
The map $\eta : \Pp \rightarrow \PP H^0(\Oo(1,2))^*$ is an isomorphism. Moreover, the restriction of $\eta$ to $Q'$ is an isomorphism to $Q$ sending an extension of type (\ref{ext3}) to $\phi(p)\in Q$.
\end{proposition}

\begin{remark}
Let $p\in \PP_3 \backslash Q$ and $\phi'$ be the restriction of the projection from $\PP_3$ to $\PP_2$ at $p$, to $Q$. For the cotangent bundle of $\PP_2$, twisted by $\Oo_{\PP_2}(2)$, admits the following exact sequence,
$$0\rightarrow \Oo_{\PP_2} \rightarrow \Omega_{\PP_2}(2) \rightarrow I_p(1) \rightarrow 0,$$
where $p$ is a point on $\PP_2$, not the point corresponding to the line passing through $p$ and the vertex point $P_0$. If we pull back the sequence via $\phi ' \circ \phi$, then we get a vector bundle $E$ admitting an exact sequence (\ref{ext4}), where $Z$ is $\phi^{-1} \circ \phi'^{-1}(p)$. This defines a map from $\PP_3$ to $\Pp$ and in fact, it extends to the inverse morphism of $\eta$.
\end{remark}

Let $C$ be a non-hyperelliptic curve of genus 4 with the two trigonal line bundles $g_3^1$ and $h_3^1$ such that $|g_3^1|=|h_3^1|$. In particular, $C$ is embedded into $\PP_3$ by the canonical embedding and there exists a unique quadric cone $Q\subset \PP_3$ containing $C$. Let $P_0$ be the vertex point of $Q$. Recall that $\FF_2$ is the minimal resolution of $Q$ at $P_0$. Let $C'$ be the proper transform of $C$ in $\FF_2$. Note that $C$ and $C'$ are isomorphic, so we will use $C$ instead of $C'$ if there is no confusion. Let us assume that the divisor type of $C\subset \FF_2$ is $(a,b)$. From the adjunction formula, we have
$$6=2g(C)-2=C.(C+K)=(a,b).(a-2,b-4).$$
Since $C$ does meet the vertex $P_0$, we have $C . \sigma=0$. Hence $(a,b)=(3,6)$.
If we tensor the following exact sequence
\begin{equation}\label{seq1}
0\rightarrow \Oo(-3,-6) \rightarrow\Oo\rightarrow \Oo_C \rightarrow 0,
\end{equation}
with a bundle $E\in \Pp$ and take the long exact sequence of cohomology, we have $h^0(E|_C)=h^0(E)=3$ since $h^1(E(-3,-6))=h^1(E)=0$.
By the adjunction formula, we have
$$\Oo_C(K_C)=\Oo(K_{\FF_2})\otimes \Oo(3,6)\otimes \Oo_C=\Oo(1,2)\otimes \Oo_C,$$
i.e. the determinant of $E|_C$ is $\Oo_C(K_C)$.
\begin{lemma}
The restriction map
$$\Phi : \Pp \dashrightarrow \Ww^2,$$
sending $E$ to $E|_C$, is well-defined.
\end{lemma}
\begin{proof}
It is enough to prove that $E|_C$ is stable. Let us assume that there exists a sub-bundle $\Oo_C(D)$ with $d=\deg (D)\geq 3$. Since the degree of $K_C-D$ is less than 4, we have $h^0(\Oo_C(D))>0$ due to the Clifford theorem \cite{H} and $h^0(E|C)=3$. Thus we can assume that $D$ is effective. Since $H^0(E)\simeq H^0(E|_C)$, $D$ can be considered as the intersection of the zero of a section of $H^0(E)$ with $C$. For a section in $H^0(E)$, let us consider an exact sequence,
$$0\rightarrow \Oo(a,b) \rightarrow E \rightarrow I_Z(1-a,2-b)\rightarrow 0,$$
where $a,b\geq 0$. From the numeric invariants of $E$ and the fact that $h^0(E)=3$, we have $(a,b)=(0,0),(1,0)$ or $(0,1)$. For a general vector bundle $E\in \Pp$, the case of $(a,b)=(0,1)$ cannot happen. Indeed, it happens only when $E\in \Dd$. Since the length of $Z$ is at most 2 in each case, $d$ must be less than 3. Hence, $E|_C$ is stable.
\end{proof}

Let $g_3^1$ be the trigonal line bundle on $C$ and we have $\Oo(0,1)|_C=g_3^1$. If $E\in \Dd$, we have an exact sequence (\ref{ext3}). If $p\not \in C$, we obtain the following exact sequence
$$0\rightarrow g_3^1 \rightarrow E|_C \rightarrow g_3^1 \rightarrow 0,$$
after tensoring with $\Oo_C$. In particular, $E|_C$ is in the same equivalent class of $g_3^1\oplus g_3^1$. If $p\in C$, then we obtain an exact sequence,
$$0\rightarrow g_3^1\otimes \Oo_C(p) \rightarrow E|_C \rightarrow g_3^1 \otimes \Oo_C(-p) \rightarrow 0,$$
which implies that $E|_C$ is not semi-stable. 
 Thus we have the following assertion.

\begin{proposition}
The restriction map $\Phi : \Pp \dashrightarrow \Ww^2$ is defined by the complete linear system $|I_C(3)|$. In particular, $\Ww^2$ is isomorphic to the \textit{Donagi-Izadi cubic threefold}.
\end{proposition}
\begin{proof}
The proof is similar with the one in \cite{Huh4}. If we choose a general hyperplane section $H\subset \Pp$, then the restriction of $\Phi$ to $H$ is not defined on 6 intersection points of $C$ with $H$. Since this indeterminacy locus lie on a conic on $H$, the blow-up of $H$ at these points is a singular cubic surface in $\PP_3$. In particular, the degree of $\Phi$ is 3.

Let $E$ be a general vector bundle in $\Ww^2$ with $h^0(E)=3$. It can be checked as in (\ref{det}) that the determinant map from $\wedge^2 H^0(E)$ to $H^0(\Oo_C(K_C))$ is injective and so we can assign a point $p_E\in \PP_3$ corresponding to the dual of the cokernel of the determinant map. This defines a map $\rho$ from $\Ww^2$ to $\PP_3$ and $\eta^{-1} \circ \rho \circ \Phi$ is the identity on $\Pp$. In particular, the dimension of $\Ww^2$ is at least 3. Conversely, the dimension of $\Ww^2$ can be shown to be at most 3 as follows: Let us assume that $E$ is the extension of $\Oo_C(K_C-D)$ by $\Oo_C(D)$, where $D$ is a divisor of $C$ with the degree $d$. Because of the stability of $E$ and the result of \cite{L}, we can assume that $d=2$ and so $h^0(\Oo_C(D))=1$. In particular, we can assume that $D$ is effective. In the extension space $\PP \Ext^1 (\Oo_C(K_C-D), \Oo_C) \simeq \PP_4$, there exists $\PP_1$-parametrization corresponding to the vector bundles $E$ with $h^0(E) \geq 3$ \cite{OPP}. Thus, we have a dominant map from a $\PP_1$-bundle over $\Sec^2(C)$ to $\Ww^2$ and so the dimension of $\Ww^2$ is at most 3. Now, we know that $h^0(I_C(3))=5$ and so $\Phi$ is given by the complete linear system $|I_C(3)|$ and the image is exactly $\Ww^2$.
\end{proof}

\bibliographystyle{amsplain}
\bibliography{F_2}
\end{document}